\documentclass{article}
\usepackage{mathtools}
\usepackage{amssymb}
\usepackage{eucal}
\usepackage{amsmath}
\begin{document}
	\begin{center}
		\Large{\bf\textbf\textbf{ Generalization of The Results on Fixed Point For Couplings on Metric Spaces }}
	\end{center}
	\noindent
	\begin{center}
		Tawseef Rashid \footnote{Department of Mathematics, Aligarh Muslim University, Aligarh-202002, India.\\       
			Email : tawseefrashid123@gmail.com} and Q. H. Khan \footnote{Department of Mathematics, Aligarh Muslim University, Aligarh-202002, India.\\ Email : qhkhan.ssitm@gmail.com}\\
	\end{center}

			\textit{The purpose of this paper is to introduce the concept of self-cyclic maps, $g$-coupling, Banach type $g$-coupling which is the generalization of couplings introduced by Choudhury et al.$\cite{bcm}$. In our main result we prove the existence theorem of coupled coincidence point for Banach type $g$-couplings which extends the results by choudhury et al. $\cite{bcm}$. 	We give an example in support of our result.}\\

	\noindent\\
	\textbf{2010MSC}: 47H10, 54H25.\\
	\begin{center}
		\textbf{1. Introduction and Preliminaries}
	\end{center} The concept of coupled fixed point was introduced in the work of Guo et al. $\cite{glk}.$ Then after the  coupled contraction mapping theorem by Bhaskar and Lakshmikantham $\cite{gbl}$, and the introduction of coupled coincidence point by Lakshmikantham.V and \'{C}iri\'{c}.L $\cite{vlc}$, the coupled fixed and coupled coincidence point results reported in papers $\cite{ha,vb,bek,bck,bcak,hjlsb,vlc,lntx,lt,mby,rsbm, whnk,wsbs,wsm}$. Kirk et al. $\cite{ksv}$ gave the concept of cyclic mappping. Recently Choudhury et al.$\cite{bcm}$ introduced the concept of couplings which are actually coupled cyclic mappings with respect to two given subsets of a metric space. Choudhury et al. $\cite{bcm}$ proved the existence of strong coupled fixed point for a coupling ( w.r.t. subsets of a complete metric space). H. Aydi et al. $\cite{haah}$ proved the theorem for existence and uniqueness of strong coupled fixed point in partial metric spaces. In this paper we extend the concept of Banach type coupling to Banach type $g$-coupling by  extending the concept of coupling to $g$-coupling. We have introduced self-cyclic maps and generalize the results by Choudhury et al $\cite{bcm}$. Now we recall some definitions.\\ 
	\noindent\\  
	\textbf{Definition $1.1$} (Coupled fixed point) $\cite{gbl}$. An element $(x,y) \in X \times X$, where $X$ is any non-empty set, is called a coupled fixed point of the mapping $F : X \times X \to X$ if $F(x,y) = x~ and~ F(y,x) =y$.\\
	\textbf{Definition $1.2$} (Strong coupled fixed point) $\cite{cm}.$   An element $(x,y) \in X \times X$, where $X$ is any non-empty set,  is called a strong coupled fixed point of the mapping $F : X \times X \to X$ if $(x,y)$ is coupled fixed point and $x = y$; that is if $F(x,x) = x$.\\
	\textbf{Definition $1.3$} (Coupled Banach Contraction Mapping) $\cite{gbl}.$ Let $(X,d)$ be a metric space. A mapping $F : X \times X \to X  $ is called coupled Banach contraction if there exists $k \in (0,1)$ s.t  $\forall ~~(x,y),(u,v) \in X \times X$, the following inequality is satisfied:
	\begin{eqnarray*}
	d(F(x,y),F(u,v)) \leq \frac{k}{2} [d(x,u) + d(y,v)].\\
	\end{eqnarray*}
	\textbf{Definition $1.4$} (Cyclic mapping)  $\cite{ksv}.$ Let $A~~ and ~~B$ be two non-empty subsets of a given set $X$. Any function $f : X \to X$ is said to be cyclic (with respect to $A$ and $B$) if 
	\begin{eqnarray*}
	f(A) \subset B ~~and~~ f(B) \subset A.
	\end{eqnarray*}
	\textbf{Definition $1.5$} (Coupling) $\cite{cm}.$ Let $(X,d)$ be a metric space and $A ~and ~B$ be two non-empty subsets of $X$. Then a function $F : X \times X \to X$ is said to be a coupling with respect to $A~and~B$ if 
	\begin{eqnarray*}
	F(x,y) \in B ~and~ F(y,x) \in A~whenever~x \in A~ and~ y \in B.
	\end{eqnarray*}
	\textbf{Definition $1.6$} (Banach type coupling) $\cite{bcm}.$ Let $A$ and $B$ be two non-empty subsets of a complete metric space $(X,d)$. A coupling $F :X \times X \to X$ is called a Banach type coupling with respect to $A$ and $B$ if it satisfies the following inequality:
	\begin{eqnarray*}
	d(F(x,y),F(u,v)) \leq \frac{k}{2} [d(x,u) +d(y,v)].
	\end{eqnarray*}
	where $x,v \in A, ~~y,u \in B ~and~k \in (0,1)$.\\
	\textbf{Definition $1.7$} (Coupled coincidence point of $F$ and $g$) $\cite{vlc}.$ An element $(x,y) \in X \times X$ is called a coupled coincidence point of the mappings $F :X \times X \to X$ and $g : X \to X$ if $F(x,y) = g(x)~and~F(y,x) = g(y)$.\\
	\textbf{Definition $1.8$} (Commutative mappings) $\cite{vlc}$. For any set $X$, we say mappings $F : X \times X \to X$ ~and ~$g :X \to X$ are commutative if
	\begin{equation*}
	g(F(x,y)) = F(g(x),g(y)), ~~\forall~~x,y \in X.
	\end{equation*}
	\begin{center}
	\textbf{2. Main Result}
	\end{center}
	Here we give some definitions in support of our main result.\\
	\textbf{Definition $2.1$} (Self-cyclic mapping). Let $A$ and $B$ be two non-empty subsets of any set $X$. Then a mapping  $g: X \to X$ is said to be self-cyclic (with respect to $A$ and $B$) if
	\begin{eqnarray*}
     g(A) \subseteq A ~and~ g(B) \subseteq B.
	\end{eqnarray*}
	\textbf{Definition $2.2$} (Strong coupled coincidence points of $F$ and $g$). A coupled coincidence point $(x,y) \in X\times X$ of $F:X \times X \to X$ and $g : X \to X$ is said to be strong coupled coincidence point of $F ~and~g$, if $x = y~ i.e.~ F(x,x) = g(x)$.\\
	\textbf{Definition $2.3$} (g-coupling).  Let $(X,d)$ be a metric space and $A$ and $B$ be two non-empty subsets of $X$. Let functions $F$ and $g$ are such that $F : X \times X \to X~ and ~g : X \to X$. Then $F$ is said to be $g$-coupling (with respect to $A$ and $B$) if 
	\begin{eqnarray*}
	F(x,y) \in g(A) \cap B   ~and ~ F(y,x) \in g(B) \cap A,  ~whenever ~ x \in A~ and~ y \in B.
	\end{eqnarray*}
	\textbf{Remark $2.4$}: It should be noted that every  $g$-coupling is a coupling but converse is not true in general.\\
	\textbf{Proof}: Let $F : X \times X \to X$ be $g$-coupling, where $g : X \to X$. Then by Definition $(2.3)$, we have
	$F(x,y) \in (g(A) \cap B)  \subset B~~and~~F(y,x) \in (g(B) \cap A)  \subset A$.\\
	this shows that $F$ is a coupling (w.r.t. $A$ ~and~$B$).\\
	Clearly converse is not true in general. If $F$ is a coupling (w.r.t $A$ ~and~ $B$), it  doesn't imply that $F$ is a $g$-coupling for any $g : X \to X$ (w.r.t. $A$~~and~~$B$) as $F(x,y)$ $\in$ $B$  doesn't imply that $F(x,y) \in g(A)$, and similarly for $F(y,x)$.\\
	\textbf{Definition $2.5$} (Banach type g-coupling). Let $A$ and $B$ be two non-empty subsets of a complete metric space $(X,d)$. Then a $g$-coupling $F: X \times X \to X$ is said to be Banach type $g$-coupling (with respect to $A$ and $B$) if the following inequality holds:
	\begin{eqnarray*}
	d(F(x,y),F(u,v)) \leq \frac{k}{2}[d(gx,gu) + d(gy,gv)].
    \end{eqnarray*}
	whenever~$ x,v \in A$,~~ $y,u \in B$ ~and $~k \in (0,1)$.
	where $g$ :$ X \to X$ is a self-cyclic mapping (with respect to $A$ and $B$).\\
	\textbf{Note :} If $g$ = $I$ (the identity mapping) which is also self-cyclic, then Banach type $g$-coupling becomes Banach type coupling.\\
	\textbf{Theorem $2.6$}: Let $A$ and $B$ be any two subsets of a complete metric space $(X,d)$. If there exists a Banach type  $g$-coupling $F : X \times X \to X$ (with respect to $A$ and$B$), where $g : X \to X$ is  self-cyclic (with respect to $A$ and~$B$). If $g(A) ~and~g(B)$ are closed subsets of $(X,d)$. Then\\
	(i) $g(A) \cap g(B) \neq \emptyset$,\\
	(ii) $F~and~g$ have a coupled coincidence point in $A \times B$, i.e. there exists $(a,b)$ $\in$ $A \times B$, s.t $F(a,b) = g(a)~ and~ F(b,a) = g(b)$.\\
	\textbf{Proof :} Here $F : X \times X \to X $ is given to be Banach type $g$-coupling, i.e.\\
	\begin{equation}
	d(F(x,y),F(u,v)) \leq \frac{k}{2}[d(gx,gu) + d(gy,gv)].\\
	\end{equation}
	where $x,v \in A,~~ y,u \in B~ and~ k \in (0,1)$.\\
	Also as $g : X \to X$ is self-cyclic (w.r.t. $A~ and~ B$), so 
	\begin{eqnarray*}
	gx,gv \in g(A) \subseteq A~ and~ gy,gu \in g(B) \subseteq B
	\end{eqnarray*}
	Let  $x_0 \in A ~and~y_0 \in B$, then by defnition of $g$-coupling, we have\\
	$F(x_0,y_0) \in g(A) \cap B ~and ~F(y_0,x_0) \in g(B) \cap A$\\
	in particular $F(x_0,y_0) \in g(A)  ~and ~F(y_0,x_0) \in g(B). $\\
	If
	\begin{eqnarray*}
	F(x_0,y_0) = g(x_0)~and~F(y_0,x_0) = g(y_0),
	\end{eqnarray*}
	  then $(x_0,y_0)$ is the coupled coincidence point of $F~and~g$, so we are done in this case.
	Otherwise $\exists~ x_1 \in A ~and~ y_1 \in B$ , s.t.
	\begin{eqnarray*}
	F(x_0,y_0) = g(x_1)~ and ~F(y_0,x_0) = g(y_1).
	\end{eqnarray*}
	Now if,
	\begin{eqnarray*}
	 F(x_1,y_1) = g(x_1) ~and~F(y_1,x_1) = g(y_1),
	\end{eqnarray*}
	 then $(x_1,y_1)$ is a coupled coincidence point of $F$ and $g$, and we are through\\
	otherwise $\exists~ x_2 ~\in A ~and ~y_2 \in B $ s.t
	\begin{eqnarray*}
	F(x_1,y_1) = g(x_2)~and~F(y_1,x_1) = g(y_2).
	\end{eqnarray*}
	continuing in this way, we get sequences $\{gx_n\}$ and $\{gy_n\}$ in $g(A)$ and $g(B)$ respectively, such that
	\begin{equation}
	g(x_{n+1}) = F(x_n,y_n)~~and~~g(y_{n+1}) = F(y_n,x_n).
	\end{equation}
	Now using $(1) ~and~ (2)$, we get
	\begin{eqnarray*}
	d(gx_1,gy_2)&=& d(F(x_0,y_0),F(y_1,x_1))\\
	 &\leq& \frac{k}{2}[d(gx_0,gy_1) + d(gy_0,gx_1)]
	\end{eqnarray*}
	and
	\begin{eqnarray*}
	d(gy_1,gx_2)&=& d(F(y_0,x_0),F(x_1,y_1))\\
	& \leq& \frac{k}{2}[d(gy_0,gx_1) + d(gx_0,gy_1)].
	\end{eqnarray*}
	from above two inequalities, we have
	\begin{eqnarray*}
	d(gx_1,gy_2) + d(gy_1,gx_2)&\leq& \frac{k}{2}[d(gx_0,gy_1) + d(gy_0,gx_1) + d(gy_0,gx_1) + d(gx_0,gy_1)]\\
	&=& k[d(gx_0,gy_1) + d(gy_0,gx_1)].
	\end{eqnarray*}
	or,
	\begin{equation}
	\frac{d(gx_1,gy_2) + d(gy_1,gx_2)}{2} \leq \frac{k}{2}[d(gx_0,gy_1) + d(gy_0,gx_1)].
	\end{equation}
	using $(1),(2)~ and~(3)$, we have
	\begin{eqnarray*}
	d(gx_2,gy_3)&=& d(F(x_1,y_1),F(y_2,x_2))\\
	&\leq& \frac{k}{2}[d(gx_1,gy_2) + d(gy_1,gx_2)]\\
	&\leq& \frac{k^2}{2}[d(gx_0,gy_1) +d(gy_0,gx_1)],
	\end{eqnarray*}
	and
	\begin{eqnarray*}
	d(gy_2,gx_3)&=& d(F(y_1,x_1),F(x_2,y_2))\\
	& \leq& \frac{k}{2}[d(gy_1,gx_2) + d(gx_1,gy_2)]\\
	& \leq&  \frac{k^2}{2}[d(gx_0,gy_1) +d(gy_0,gx_1)].
	\end{eqnarray*}
	Let for some integer n,
	\begin{equation}
	d(gx_n,gy_{n+1}) \leq \frac{k^n}{2}[d(gx_0,gy_1) +d(gy_0,gx_1)].
	\end{equation}
	\begin{equation}
	d(gy_n,gx_{n+1}) \leq \frac{k^n}{2}[d(gx_0,gy_1) +d(gy_0,gx_1)].
	\end{equation}
	Now using $(1),(4)~and~(5)$, we have
	\begin{eqnarray*}
	d(gx_{n+1},gy_{n+2}) &=& d(F(x_n,y_n),F(y_{n+1},x_{n+1}))\\
	& \leq& \frac{k}{2}[d(gx_n,gy_{n+1}) + d(gy_n,gx_{n+1})]\\
	& \leq& \frac{k}{2}{k^n}[d(gx_0,gy_1) +d(gy_0,gx_1)]\\
	& = &\frac{k^{n+1}}{2}[d(gx_0,gy_1) +d(gy_0,gx_1)].
	\end{eqnarray*}
	similarly as above  we can show that 
	\begin{eqnarray*}
	d(gy_{n+1},gx_{n+2}) \leq \frac{k^{n+1}}{2}[d(gx_0,gy_1) +d(gy_0,gx_1)].
	\end{eqnarray*}
	Thus $(4)~and~(5)$ remains also true for $n+1$,\\
	Hence by principle of mathematical induction, we have $\forall~ n \geq 1$,
	\begin{equation}
		d(gx_n,gy_{n+1}) \leq \frac{k^n}{2}[d(gx_0,gy_1) +d(gy_0,gx_1)].
	\end{equation}
	\begin{equation}
	d(gy_n,gx_{n+1}) \leq \frac{k^n}{2}[d(gx_0,gy_1) +d(gy_0,gx_1)].
	\end{equation}
	Again by $(1) ~and~(2)$, we have
	\begin{eqnarray*}
	d(gx_1,gy_1) &=& d(F(x_0,y_0),F(y_0,x_0))\\
	& \leq& \frac{k}{2}[d(gx_0,gy_0) + d(gy_0,gx_0)]\\
	& =& kd(gx_0,gy_0).
	\end{eqnarray*}
	that is,
	\begin{equation}
	d(gx_1,gy_1) \leq kd(gx_0,gy_0).
	\end{equation}
	Then from $(1),(2)~and~(8)$, we have
	\begin{eqnarray*}
	d(gx_2,gy_2) &=& d(F(x_1,y_1),F(y_1,x_1))\\
	& \leq &\frac{k}{2}[d(gx_1,gy_1) + d(gy_1,gx_1)]\\
	& = &kd(gx_1,gy_1)\\
	& \leq& k^2d(gx_0,gy_0).
	\end{eqnarray*}
	Let for some integer n, we have
	\begin{equation}
	d(gx_n,gy_n) \leq k^nd(gx_0,gy_0).
	\end{equation}
	Then from $(1),(2)~and~(9)$, we get
	\begin{eqnarray*}
	d(gx_{n+1},gy_{n+1}) &=& d(F(x_n,y_n),F(y_n,x_n))\\
	& \leq& \frac{k}{2}[d(gx_n,gy_n) + d(gy_n,gx_n)]\\
	& = &kd(gx_n,gy_n)\\
	& \leq &k^{n+1}d(gx_0,gy_0).
	\end{eqnarray*}
	This shows that $(9)$ remains also true for $n+1$, thus by principle of mathemtical induction, we say that\noindent
	\begin{equation}
	d(gx_n,gy_n) \leq k^nd(gx_0,gy_0),~ \forall~n \geq 1.
	\end{equation}
	Now by $(6),(7)~and~(10)$ and by triangular inequality, we have $\forall~n \geq 1$
	\begin{eqnarray*}
	d(gx_n,gx_{n+1}) + d(gy_n,gy_{n+1}) &\leq& d(gx_n,gy_n) + d(gy_n,gx_{n+1}) + d(gy_n,gx_n) + d(gx_n,gy_{n+1})\\
	& = &2d(gx_n,gy_n) + [d(gy_n,gx_{n+1}) + d(gx_n,gy_{n+1})]\\
	& \leq& 2k^nd(gx_0,gy_0) + k^n[d(gx_0,gy_1) + d(gy_0,gx_1)].
	\end{eqnarray*}
	since $k \in (0,1)$, it follows that 
	$\sum d(gx_n,gx_{n+1}) + \sum d(gy_n,gy_{n+1}) < \infty$.\\
	Thus sequences $\{gx_n\}$ and $\{gy_n\}$ are Cauchy sequences in $g(A)$ and $g(B)$.\\
	As $g(A)$ and $g(B)$ are closed subsets of complete metric space $(X,d)$, so $\{g(x_n)\}~and~\{g(y_n)\}$ are convergent in $g(A) ~and ~g(B)$ respectively.\\
	Therefore $\exists~u \in g(A) ~and~v \in g(B), s.t. $
	\begin{equation}
	gx_n \to u ~and ~gy_n \to v ~as~ n \to \infty.
	\end{equation}
	using $(10)$, as $k \in (0,1)$, we have
	\begin{eqnarray*}
	d(gx_n,gy_n) \to 0 ~~~~as~ n \to \infty.
	\end{eqnarray*}
	therefore from $(11)$, we have
	\begin{eqnarray}
	u = v.
	\end{eqnarray}
	As $u \in g(A) ~and ~v \in g(B) \Rightarrow u \in g(A) \cap g(B)$\\
	This proves part (i) that $g(A) \cap g(B) \neq \emptyset$.\\
	Now, since $u \in g(A) ~and~v \in g(B)$, therefore $\exists~~ a \in A~~and~~b \in B$, s.t\\
	$u = g(a)~and~v = g(b)$.\\
	then from $(11) ~and~(12)$, we have
	\begin{eqnarray}
	gx_n \to g(a) ~~and~~ gy_n \to g(b).
	\end{eqnarray}
	and
	\begin{equation}
	g(a) = g(b).
	\end{equation}
	Now by $(1),(2),(13),(14)$ ~and~triangular inequality, we have
	\begin{eqnarray*}
	d(g(a),F(a,b)) &\leq &d(g(a),gy_{n+1}) + d(gy_{n+1},F(a,b))\\
	& =& d(g(a),gy_{n+1}) + d(F(y_n,x_n),F(a,b))\\
	& \leq& d(g(a),gy_{n+1}) + \frac{k}{2}[d(gy_n,g(a)) + d(gx_n,g(b))]\\
	& \to &0 ~~~~as~~~~n \to \infty.
	\end{eqnarray*}
	thus, we have
	\begin{equation}
	F(a,b) = g(a).
	\end{equation}
	Again from$(1),(2),(13),(14)$ and traingular inequality, we have
	\begin{eqnarray*}
	d(g(b),F(b,a)) &\leq &d(g(b),gx_{n+1}) + d(gx_{n+1},F(b,a))\\
	&=& d(g(b),gx_{n+1}) + d(F(x_n,y_n),F(b,a))\\
	& \leq &d(g(b),gx_{n+1}) + \frac{k}{2}[d(gx_n,g(b)) + d(gy_n,g(a))]\\
	&\to& 0 ~~~~as~~~~n \to \infty. 
	\end{eqnarray*}
	thus, we get
	\begin{equation}
	F(b,a) = g(b).
	\end{equation}
	Hence from $(15) ~and~(16)$ we get,
	$F(a,b) = g(a)$~~and~~$F(b,a) = g(b)$,\\  where $a \in A~~and~~b \in B$.\\
	Thus $(a,b) \in A \times B$ is the coupled coincidence point of $F~and~ g$.\\
	\textbf{Remark $2.7$}. It is worth noting that the above theorem also gives the condition for the existence of symmetric point of $F ~ in ~A \times B$, i.e. $\exists ~(a,b)\in A \times B$  s.t \\
	$F(a,b) = F(b,a)$.\\ As from $(14)$~ $g(a) = g(b)$~  so,  $F(a,b) = F(b,a)$.\\
	\noindent\\ 
	\textbf{Existence and Uniqueness of Strong coupled coincidence point of $F$ and g}.\\
	\textbf{Theorem $2.8$}: If in addition to above condition in Theorem $2.6$, $g$ is one-one, then \\
	(i) $A \cap B \neq \emptyset$, and\\
	(ii) $F ~and~g$ have unique strong coupled coincidence point in $A \cap B$.\\
	\textbf{Proof:} since $g$ is given to be one-one\\
	then from $(14)$ of Theorem $2.6$, we get
	\begin{eqnarray*}
	a = b.
	\end{eqnarray*} 
	since, $a \in A ~~and~~ b \in B$\\
	$\Rightarrow a = b \in A \cap B$.
	Hence $A \cap B \neq \emptyset$, this proves our part (i).\\
	Also from $(15)$ of Theorem $2.6$, we have
	\begin{equation*}
	F(a,a) =  g(a).
	\end{equation*}
	which shows that $F$ and $g$ have strong coupled coincidence point in $A \cap B$.\\
	\textbf{Uniqueness:} Let us suppose if possible there exists two strong coupled coincidence points $l,m \in A \cap B$ of $F~~and~~g$\\
	then by definition we have,
	\begin{equation}
	F(l,l) = g(l)~~and~~ F(m,m) = g(m).
	\end{equation}
	then from $(1)$ of Theorem $2.6$, we have
	\begin{eqnarray*}
	d(g(l),g(m)) &= &d(F(l,l),F(m,m))\\
	& \leq &\frac{k}{2}[d(g(l),g(m)) + d(g(l),g(m))]\\
	& =& k[d(g(l),g(m))].	
	\end{eqnarray*}
	Which is a contradiction as $k \in (0,1)$ and is only possible if $d(g(l),g(m)) = 0$, 
	i.e. $g(l) = g(m)$ 
	or $l = m$ because $g$ is one-one.\\
	Hence $F$ and $g$ have a unique strong coupled coincidence point in $A \cap B$.\\
	\textbf{Corrollary $2.9$ :} If $F : X \times X \to X$ and $g : X \to X$ are commutative and $g$ is one-one.\\
	if $(gx,gy)$ is a coupled coincidence point of $F~and ~g$, where $x,y \in X$, then\\
	$(x,y)$ is also a coupled coincidence point of $F~~and~~g$.\\
	\textbf{Proof :} As $(gx,gy)$ is a coupled coincidence point of $F ~and~g$, we have
	\begin{eqnarray}
       F(gx,gy) = gx,\\
	   F(gy,gx) = gy.
	\end{eqnarray}
	Now by ($18$)~and commutavity of $F~and ~g$, we have
	\begin{eqnarray*}
	g(F(x,y)) &=& F(gx,gy)\\
	& =& gx.
	\end{eqnarray*}
	as $g$ is one-one, we have from above
	\begin{equation}
	F(x,y) = x.
	\end{equation}
		Now again by ($19$)~and commutavity of $F~and ~g$, we have\\
		\begin{eqnarray*}
			g(F(y,x)) &=& F(gy,gx)\\
			& = &gy.
		\end{eqnarray*}
		as $g$ is one-one, we have from above
		\begin{equation}
		F(y,x) = y.
		\end{equation}
		Thus from $(20)~and~(21)$, we get
		\begin{eqnarray*}
		F(x,y) = x ~~and~~F(y,x) = y.
		\end{eqnarray*}
		i.e. $(x,y)$ is the coupled coincidence point of $F ~~and~~g$.\\
		The following example illustrates our results.\\
		\textbf{Example $2.10$}. Let $X = R$ with the metric defined as $d(x,y)  = \mid x - y\mid$, where $x,y \in X$.\\
		Let $A = [0,2] ~and~B = [0,3]$.\\
		Let $F$ be defined as $F(x,y) = \frac{x+y}{10}$, where $x,y \in X$.\\
		and let $g : X \to X$ is defined by $g(x) = \frac{x}{2}$.\\
		Then $g(A) = [0,1] ~~and~~g(B) = [0,\frac{3}{2}]$, so $g(A) ~and~g(B)$ are closed subsets of $X$.\\
		Also $g(A) \subseteq A~~and ~~g(B) \subseteq B$, so $g$ is self-cyclic.\\
		Now we show that $F$ is $g$-coupling.\\
		As $g(A) \cap B = [0,1]$ ~~and ~~$g(B) \cap A = [0,1]$, so  $\forall~ x \in A ~~and~~y \in B$, we have\\
		$0 \leq F(x,y) \leq \frac{1}{2}$~~and~~$0 \leq F(y,x) \leq \frac{1}{2}$.\\
		i.e. $F(x,y) \in g(A) \cap B$ ~~and~~$F(y,x) \in g(B) \cap A$, which shows that $F$ is $g$-coupling with respect to $A~~and~~B$.\\
		Again for $x,v \in A~~and~~y,u \in B$, we have
		\begin{eqnarray*}
		d(F(x,y),F(u,v)) &=& \mid \frac{x + y}{10} - \frac{u + v}{10} \mid\\	
		&=& \frac{1}{5} \mid \frac{x-u}{2} + \frac{y-v}{2} \mid \\
		& \leq& \frac{1}{5}[\mid \frac{x-u}{2} \mid  + \mid \frac{y-v}{2} \mid]\\
		& =& \frac{1}{5}[\mid \frac{x}{2} - \frac{u}{2} \mid  +  \mid \frac{y}{2} - \frac{v}{2} \mid ]\\
		& =& \frac{k}{2}[d(gx,gu) + d(gy,gv)],where~~ k = \frac{2}{5} \in (0,1).
		\end{eqnarray*}
		which shows that $F$ is a Banach type $g$-coupling (w.r.t. A and B).\\
		Thus all the conditions of Theorem $(2.6)$ are satisfied, therefore $\exists~~ (a,b) \in A \times B$ s.t.\\
		$F(a,b) = g(a)~~and~~F(b,a) = g(b)$.\\
		i.e. 
		\begin{equation}
		 \frac{a+b}{10} = \frac{a}{2}~~and~~\frac{b+a}{10} = \frac{b}{2}
		\end{equation}
		or
		\begin{eqnarray*}
		\frac{a+b}{5} = a = b.
		\end{eqnarray*}
		which is only possible with $a =0~~and~~b = 0$.\\
		thus $(0,0)$ is the unique strong $g$-coupled coincidence point of $F$ and $g$.\\
		\textbf{Note :} The uniqueness of the strong $g$-coupled coincidence point in the above example is because of $g$ is  one-one.\\
		\textbf{Corrollary $2.11$:} Let $g :X \to X $ is a Banach  contraction, Then every Banach type $g$-coupling (w.r.t. A~and~B) is a Banach type coupling (w.r.t. A and B), where $A~and~B$ are subsets of $X$.\\
		\textbf{Proof :} since $g$ is a contraction, therefore $ \exists~ \alpha \in (0,1)$, s.t
		\begin{equation}
		d(gx,gy) \leq \alpha d(x,y)~~~ \forall~ x,y \in X.
		\end{equation}
		Let $F :X \times X \to X$ be a Banach type $g$-coupling (w.r.t. $A$ and $B$), i.e.
		\begin{equation}
		d(F(x,y),F(u,v)) \leq \frac{k}{2}[d(gx,gu) + d(gy,gu)],~ where~ x,v \in A, ~ ~y,u \in B~ and~ k \in (0,1). 
		\end{equation}
		Also $F$ is a $g$-coupling (w.r.t. $A$ and $B$) and hence by remark ($2.4$), $F$ is a coupling (w.r.t. $A$ and $B$).\\
		Now using $(23)$ in $(24)$, we get
		\begin{eqnarray*}
		d(F(x,y),F(u,v)) &\leq &\frac{k}{2}[\alpha d(x,u) + \alpha d(y,v)]\\
		& \leq &\frac{\alpha k}{2}[d(x,u) + d(y,v)].
		\end{eqnarray*}
		as $k \in (0,1) ~and~ \alpha \in (0,1), ~therefore~~ k\alpha  = k_1\in (0,1)$, so
		\begin{eqnarray*}
		d(F(x,y),F(u,v)) \leq \frac{k_1}{2}[d(x,u) + d(y,v)].
		\end{eqnarray*}
		where $x,v \in A, ~~ y,u \in B ~~and~~ k_1 \in (0,1)$.\\
		Hence $F$ is a Banach type coupling (w.r.t. A and B).
		\begin{center}
			\textbf{3. Conclusion}
		\end{center}
		 In this paper we introduce some new results that extend the concept of coupling and Banach type coupling  introduced by Choudhury et al.  $\cite{bcm}$ to $g$-coupling and Banach type $g$-coupling resp. and proved the existence and uniqueness theorem for coupled coincidence point and strong coupled coincidence point.\\

\end{document}